\documentclass[11pt]{amsart}

\usepackage{amsmath,amsthm,amsfonts,amssymb,bm,graphicx}

\usepackage
{hyperref}
\hypersetup{colorlinks=true,citecolor=blue,linkcolor=blue,urlcolor=blue,
pdfstartview=FitH, pdfauthor=Valentin Blomer ,pdftitle=Period integrals and Rankin-Selberg $L$-functions on $GL(n)$}

\theoremstyle{plain}
\newtheorem{theorem}{Theorem}
\newtheorem{cor}{Corollary}
\newtheorem{lemma}{Lemma}

\numberwithin{equation}{section}

\theoremstyle{definition}

\renewcommand{\geq}{\geqslant}
\renewcommand{\leq}{\leqslant}

\makeatletter
\DeclareRobustCommand\widecheck[1]{{\mathpalette\@widecheck{#1}}}
\def\@widecheck#1#2{%
    \setbox\z@\hbox{\m@th$#1#2$}%
    \setbox\tw@\hbox{\m@th$#1%
       \widehat{%
          \vrule\@width\z@\@height\ht\z@
          \vrule\@height\z@\@width\wd\z@}$}%
    \dp\tw@-\ht\z@
    \@tempdima\ht\z@ \advance\@tempdima2\ht\tw@ \divide\@tempdima\thr@@
    \setbox\tw@\hbox{%
       \raise\@tempdima\hbox{\scalebox{1}[-1]{\lower\@tempdima\box
\tw@}}}%
    {\ooalign{\box\tw@ \cr \box\z@}}}
\makeatother

\begin{document}

\author{Valentin Blomer}
\address{Mathematisches Institut, Bunsenstr. 3-5, 37073 G\"ottingen, Germany} \email{blomer@uni-math.gwdg.de}

\title{Period integrals and Rankin-Selberg $L$-functions on $GL(n)$}
 
\thanks{The author was supported   by a Volkswagen Lichtenberg Fellowship and a Starting Grant of the European Research Council.  }


\begin{abstract} We compute the second moment of a certain family of Rankin-Selberg $L$-functions $L(f \times g, 1/2)$ where $f$ and $g$ are Hecke-Maass cusp forms on $GL(n)$. Our bound is as strong as the Lindel\"of hypothesis on average, and recovers individually the convexity bound. This result is new even in the classical case $n=2$.
\end{abstract}


\setcounter{tocdepth}{2}  \maketitle 

\section{Introduction}

Automorphic $L$-functions are naturally given as Dirichlet series and can therefore be investigated using the arithmetic properties of their coefficients. Often, however, $L$-functions  can be interpreted as period integrals over the associated automorphic forms;  this gives a more geometric approach to $L$-functions and moments thereof. Spectacular successes by such techniques have recently been obtained, for example, in \cite{V, BR, MV}, yielding  very strong and general subconvexity bounds for $L$-functions of degree 2, 4 and 8. $L$-functions and automorphic forms for higher rank groups are quite mysterious objects, and their analytic properties are not well understood in many respects. This paper  establishes sharp (Lindel\"of-type) bounds for second moments of certain families of $L$-functions of arbitrarily high degree. \\

The last ten years have seen  a large number of deep and technically involved results on classical Rankin-Selberg $L$-functions in various aspects, starting with \cite{Sa, KMV} and culminating in  \cite{HaM, LLY, JM, MV}. Typically one of the factors is fixed, while a subconvex bound is obtained with respect to one or  more parameters of the other factor. Here we are interested in the situation when both factors have varying parameters (cf.\ the preprint \cite{HoM} for a very different scenario of this type). We start with the most classical case: holomorphic cuspidal Hecke eigenforms $f, g \in S_k$ of even weight $k$ for the group $\Gamma = SL_2(\Bbb{Z})$. 

\begin{theorem} Let $f \in S_k$ be  Hecke eigenform, and let $B_k$ be an orthogonal Hecke basis of $S_k$. Then  
\begin{displaymath}
  \sum_{g \in B_k} |L(f \times \bar{g}, 1/2)|^2 \ll_{\varepsilon} k^{1+\varepsilon}
\end{displaymath}
for any $\varepsilon >0$. 
\end{theorem}

Note that $\text{dim}_{\Bbb{C}} S_k \sim k/12$, so this bound is as strong as the Lindel\"of conjecture on average. Dropping all but one term, we recover the convexity bound $L(f \times \bar{f}, 1/2) \ll k^{1/2+\varepsilon}$ and also $L(\text{sym}^2 f, 1/2) \ll k^{1/2+\varepsilon}$. Moreover, the Lindel\"of hypothesis is true for almost all $g$ in the following sense:

\begin{cor} Let  $f \in S_k$ be  Hecke eigenform and let $\delta > 0$. Then the bound
\begin{displaymath}
  |L(f \times \bar{g}, 1/2)| \leq k^{\delta}
\end{displaymath}
is satisfied for all but $O(k^{1-2\delta+\varepsilon})$ Hecke eigenforms $g \in B_k$.
\end{cor}

A classical proof of Theorem 1 would start with an approximate functional equation, followed by an application of Petersson's trace formula.  This leads to terms of the form
\begin{displaymath}
   \sum_{m, n \asymp k} \lambda_f(n) \bar{\lambda}_f(m)   J_{k-1}\left( 4\pi \sqrt{nm}\right)
\end{displaymath}
and slightly more complicated expressions. 
Note that we are facing the  complicated asymptotic behaviour of the Bessel function in the transitional range, and that the weight of the Fourier coefficients involved in this sum is large. Nevertheless, one would still try to apply Voronoi summation in one of the variables which   introduces another $J_{k-1}$-function. One would then hope for a miracle (some special formula for products of Bessel functions) to obtain the bound $k^{1+\varepsilon}$.

Perhaps this argument can be pushed through, but we will present in the next section a very short and clean proof of Theorem 1 based on a period integral approach that completely avoids trace formulae, approximate functional equations, Voronoi  summation etc. The same argument works for Maa{\ss} forms of large spectral parameter: if $f$ is a Maa{\ss} form for $SL_2(\Bbb{Z})$ with large spectral parameter $\nu \in i \Bbb{R}$, then
\begin{equation}\label{Maass}
  \sum_{|\mu - \nu | \leq 1} |L(f \times \bar{g}, 1/2)|^2 \ll |\nu|^{1+\varepsilon}
\end{equation}
where the sum is over an orthogonal basis of  Hecke-Maa{\ss} cusp forms $g$ for $SL_2(\Bbb{Z})$ with spectral parameter $\mu$ satisfying $|\mu - \nu| \leq 1$. In fact, the method is sufficiently strong and flexible to generalize to $GL_n$ for arbitrary $n$. Here explicit trace formulae and such tools are not even available, and a conventional approach would be hopeless. We proceed to describe the general result in more detail.\\
 
Let $n \geq 2$ and let $f$ be a tempered, spherical Hecke-Maa{\ss} form for the group $SL_n(\Bbb{Z})$ with spectral parameters $\nu = (\nu_1, \ldots, \nu_{n-1}) \in i\Bbb{R}^{n-1}$. (We follow the notation of \cite{Go}, except that  our unitary axis is $i\Bbb{R}$ rather than $1/n + i \Bbb{R}$.)  Let $g$ be another tempered Maa{\ss}  for $SL_n(\Bbb{Z})$ whose spectral parameters we generally denote by $\mu = (\mu_1, \ldots, \mu_{n-1}) \in i\Bbb{R}^{n-1}$. 
 Let 
\begin{displaymath}
  \Gamma_{\Bbb{R}}(s) = \pi^{-s/2} \Gamma(s/2). 
\end{displaymath}
We define the measure
\begin{equation}\label{defspec}
  d_{\text{spec}}\mu = \prod_{1 \leq j \leq k \leq n-1} G\left(n(\mu_j + \ldots + \mu_k)\right) d\mu, \quad G(ix) := \left|\frac{\Gamma_{\Bbb{R}}(1 +  ix)}{\Gamma_{\Bbb{R}}( i x)}\right|^2 = \frac{x}{2\pi} \tanh\left(\frac{\pi x}{2}\right).
  \end{equation}
Up to a positive constant, this is the Plancherel measure on $SL_n$, measuring the density of Maa{\ss} forms (see Section 3). 
For example, if $n=3$ then $d_{\text{spec}}\mu$ is roughly $|\mu_1\mu_2(\mu_1+\mu_2)|  d\mu$. 

\begin{theorem} We have
\begin{displaymath}
  \sum_{\| \mu - \nu \| \leq 1} |L(f \times \bar{g}, 1/2)|^2 \ll \Bigl(\int_{\| \mu - \nu \| \leq 1}   d_{\text{spec}} \mu\Bigr)^{1+\varepsilon}
\end{displaymath}
where the sum runs over an orthogonal basis of tempered Hecke-Maa{\ss} cusp forms $g$ for the group $SL_n(\Bbb{Z})$ with spectral parameter $\mu \in i \Bbb{R}^{n-1}$ satisfying $\|\mu - \nu \| \leq 1$. 
\end{theorem}


The temperedness assumption  is only made for convenience and not essential. Again Theorem 2 is as strong as the Lindel\"of hypothesis on average, and it is at the edge of subconvexity: dropping all but one term, one recovers the convexity bound for $L(f \times \bar{f}, 1/2)$.\footnote{We remark, however, that from a logical point of view Theorem 2 does not give a new proof of the convexity bound for Rankin-Selberg $L$-functions, as the convexity bound in the form of  \cite[Theorem 2]{Li}  is implicitly used in the argument. Nevertheless, this remark sheds light on the quality of the mean value result in Theorem 2.} It also shows that for \emph{almost all} Maa{\ss} forms $g$ for  $SL_n(\Bbb{Z})$ with $\|\mu - \nu\| \leq 1$ the central value $L(f \times \bar{g}, 1/2)$ satisfies the Lindel\"of hypothesis. 
\\
 

Theorem 2  is the first time that moments of Rankin-Selberg $L$-functions or any other type of  $L$-functions in   arbitrary  rank are estimated efficiently. 
The key of success is the geometric approach to $L$-functions via periods of automorphic forms that never touches its Fourier coefficients explicitly. 

It will be clear from the proof that the bounds of Theorems 1 and 2 hold for any (fixed) point $s=1/2+it$ on the critical line. Moreover, at least in the situation of Theorem 1 the factor $k^{\varepsilon}$ can be replaced by some power of $\log k$. Another potential situation for the application of the method is the level aspect. In the situation of Theorem 1 it can likely provide bounds of the type $\sum_{g \in B_k(q)} |L(f \times \bar{g}, 1/2)|^2 \ll_k q^{1+\varepsilon}$ where $f$ and $g$ are now of (squarefree) level $q$. We leave this and other applications to future work. \\

In rank 1 there are more period formulas available, and we remark that the following variant of \eqref{Maass} holds. 

\begin{theorem} Let $h$ be a fixed Hecke-Maa{\ss} cusp form for $SL_2(\Bbb{Z})$ and $f$ a Hecke-Maa{\ss} cusp form with large spectral parameter $\nu \in i \Bbb{R}$, then  one has the following  Lindel\"of-on-average estimate for triple product $L$-functions:
\begin{displaymath}
  \sum_{|\mu - \nu | \leq 1} L(f \times \bar{g} \times h, 1/2) \ll_h |\nu|^{1+\varepsilon}
\end{displaymath}
where the sum runs over all $g$ satisfying the same summation condition as in \eqref{Maass}. 
\end{theorem} 

Note that the central value is non-negative.  Dropping all but the term $g=f$, this recovers the convexity bound for $L(f \times \bar{f} \times h, 1/2)$. (The work of Bernstein-Reznikov establishes subonvexity if only \emph{one} of the three factors in $L(f \times g \times h, 1/2)$ has large spectral parameter.) 


\section{Proof of Theorem 1}

Let  
\begin{equation}\label{four}
  f(z) =   \sum_{n=1}^{\infty} \lambda_f(n) (4 \pi n)^{(k-1)/2} e(nz) \Gamma(k)^{-1/2} \in S_k, \quad z \in \Bbb{H},
\end{equation}
be a cuspidal holomorphic Hecke eigenform of even weight $k$ for the group $\Gamma = SL_2(\Bbb{Z})$ with Hecke eigenvalues $\lambda_f(n)$. For convenience we have included a normalizing Gamma factor\footnote{The double use of the symbol $\Gamma$ as the Gamma-function and the modular group should not lead to confusion.}. The space $S_k$ is a finite dimensional Hilbert space (of dimension $\sim k/12$) with inner product
 \begin{displaymath}
   \langle f, g\rangle = \int_{\Gamma \backslash \Bbb{H}} f(z) \bar{g}(z) y^k \frac{dxdy}{y^2}.
 \end{displaymath}
For $z \in \Bbb{H}$ let 
\begin{displaymath}
  E(z, s) :=  \sum_{\gamma \in \overline{\Gamma}_{\infty} \backslash \overline{\Gamma}} \Im (\gamma z)^s =  \frac{1}{2} \sum_{\gamma \in \Gamma_{\infty} \backslash \Gamma} \Im (\gamma z)^s
\end{displaymath}
be the standard Eisenstein series, initially defined in $\Re s > 1$, with meromorphic continuation to all $s \in \Bbb{C}$. (Here $\overline{\Gamma}$ is the image of $\Gamma$ in $PSL_2(\Bbb{R})$ and $\Gamma_{\infty}$ is the group of integral unipotent upper triangular matrices.)  It has a simple pole at $s=1$ with constant residue $3/\pi$. By the Rankin-Selberg unfolding method we can compute the norm of $f$:
\begin{displaymath}
\begin{split}
 \| f\|^2 & = \frac{\pi}{3} \underset{s=1}{\text{res}} \int_{\Gamma\backslash \Bbb{H}} |f(z)|^2 E(z, s) \frac{dx dy}{y^2}  =   \frac{\pi}{3 \Gamma(k)} \underset{s=1}{\text{res}} \int_0^{\infty}  \sum_{n=1}^{\infty} |\lambda_f(n)|^2 (4\pi n)^{k-1} e^{-4\pi n y}  y^{s+k} \frac{dy}{y^2} \\
 &= \frac{\pi}{3\Gamma(k)} \underset{s=1}{\text{res}} \frac{1}{\zeta(2s)} L(f\times \bar{f}, s) \frac{\Gamma(s+k-1)}{(4\pi)^s} = \frac{1}{2\pi^2}  \underset{s=1}{\text{res}} L(f \times \bar{f}, s)
\end{split} 
\end{displaymath}
where $L(f \times \bar{f}, s) = \zeta(2s) \sum_n |\lambda_f(n)|^2 n^{-s}$ is the Rankin-Selberg $L$-function. It well-known that $L(f \times \bar{f}, s)$ in $\Re s \geq 1+\varepsilon$ as well as the residue at $s=1$ are uniformly bounded by $k^{\varepsilon}$. This follows either by using Deligne's bounds for $\lambda_f(n)$ or more elementarily  from a trick of Iwaniec that is described in \cite[p.\ 119-120]{Iw} for non-holomorphic cusp forms.\\

Let $E^{\ast}(z, s) := \Gamma_{\Bbb{R}}(2s)\zeta(2s) E(z, s)$. By Bessel's inequality\footnote{We may not have an equality, because the spectrum also contains weight $k$ Maa{\ss} forms and Eisenstein series; I thank G.\ Harcos for pointing this out.} we have
\begin{displaymath}
\begin{split}
  \| f E^{\ast} (., s) \|^2 & \geq \sum_{g \in B_k} \frac{1}{\| g\|^2} |\langle f E^{\ast}(., s), g \rangle|^2 
\end{split}  
\end{displaymath}
for all $s \not= 1$. Unfolding yields (first in $\Re s > 1$, but then everywhere by analytic continuation)
\begin{displaymath}
  \langle f E^{\ast}(., s), g \rangle =    L(f \times \bar{g}, s) \frac{\Gamma(s+k-1)\Gamma(s)}{(2\pi)^{2s}\Gamma(k)} .
\end{displaymath}
 Specializing to $s = 1/2$, we obtain
\begin{displaymath}
 \| f E^{\ast}(., 1/2)\|^2 \gg    \sum_{g \in B_k}  \frac{1}{\| g\|^2} \Bigl| L(f \times \bar{g}, 1/2) \frac{\Gamma( k-\frac{1}{2}) }{\Gamma(k)} \Bigr|^2   \gg   \frac{1}{k^{1+\varepsilon}} \sum_{g \in B_k} |L(f \times \bar{g}, 1/2)|^2, 
\end{displaymath}
and it remains to bound the left hand side. From the Fourier expansion of $E^{\ast}(z, 1/2)$ it is easy to see (cf.\ e.g.\ \cite[p.\ 61]{Iw}) that $E^{\ast}(z, 1/2) \ll y^{1/2}(1+ |\log y|) \ll y^{1/2 +\varepsilon}$ for $y\geq 1/2$. Let $\mathcal{F}$ denote the standard fundamental domain for $\Gamma \backslash \Bbb{H}$. Then
\begin{equation}\label{trick}
 \|  f E^{\ast}(., 1/2)\|^2 \ll \int_{\mathcal{F}} |f(z)|^2  y^{1+\varepsilon}  y^{k} \frac{dx dy}{y^2} \leq \int_{\mathcal{F}} |f(z)|^2  E(z, 1+\varepsilon)  y^{k} \frac{dx dy}{y^2}.
 \end{equation}
 Unfolding once again, the right hand side equals
 \begin{displaymath}
  \sum_{n=1}^{\infty} \frac{|\lambda_f(n)|^2}{n^{1+\varepsilon}} \frac{\Gamma(k+\varepsilon)}{(4\pi)^{1+\varepsilon}\Gamma(k)} \ll k^{\varepsilon}, 
 \end{displaymath}
 and the proof is complete.\\
 
\emph{Remark:} Notice how $|E(., 1/2)|^2$ is transformed into $E(., 1+\varepsilon)$ in \eqref{trick}. This feature is also apparent in \cite{MV} in the course of regularizing $E(., 1/2)$. We observe that the  proof of Theorem 1 is captured in  the chain
\begin{equation}\label{backbone}
\begin{split}
 & \frac{1}{k^{1+\varepsilon}} \sum_{g \in B_k} |L(f \times \bar{g}, 1/2)|^2 \ll  \sum_{g \in B_k} \frac{1}{\|g\|^2} |\Lambda^{\ast}(f \times \bar{g}, 1/2)|^2 = \sum_{g \in B_k} \frac{1}{\|g\|^2} |\langle f E^{\ast}(., 1/2), g\rangle|^2 \\
 & \leq   \| f E^{\ast}(., 1/2) \|^2 \ll \langle f E^{\ast}(., 1+\varepsilon), f\rangle =  \Lambda^{\ast}(f\times \bar{f}, 1+\varepsilon) \ll k^{\varepsilon}
\end{split}  
\end{equation} 
where 
\begin{displaymath}
  \Lambda^{\ast}(f \times \bar{g}, s) = \frac{\Lambda(f\times \bar{g}, s)}{ \Gamma(k)} = L(f \times \bar{g}, s) \frac{\Gamma(s+k-1)\Gamma(s)}{(2\pi)^{2s} \Gamma(k)}
\end{displaymath}  
   denotes a  ``re-normalized" completed $L$-function (so that $\| f \|, \| g \| \approx 1$). We shall see that the same argument works for arbitrary rank, the only ``hard" ingredients being Li's \cite{Li}  result on uniform bounds for Rankin-Selberg $L$-functions close to $s=1$, and Stade's formula that is needed twice for the first inequality and once for the last inequality in \eqref{backbone}. The regularization of the $GL_n$ Eisenstein series is carried out in Lemma 1 below which is of independent interest and  may have applications in other situations. \\
   

\section{Automorphic forms on $GL(n)$}  

Let $SL_n(\Bbb{R}) = NAK$ be the Iwasawa decomposition, let $W$ be the Weyl group, and let $\mathfrak{a}$ be the Lie algebra of $A$. We can view a tempered Maa{\ss} form $f$ for the group $SL_n(\Bbb{Z})$ as an element of $i\mathfrak{a}^{\ast}/W$; the corresponding linear form $l = (\alpha_1, \ldots, \alpha_{n-1})$   contains the $n-1$ archimedean Langlands parameters of $f$. A convenient basis in $\mathfrak{a}_{\Bbb{C}}^{\ast}$ is given by the fundamental weights; the coefficients of $l$ with respect to this basis can be obtainted by evaluating $l$ at the co-roots $\text{diag}(1, -1, 0, \ldots, 0), \ldots, \text{diag}(0, \ldots, 0, 1, -1)$, giving $n-1$ numbers that we call $n\nu_1, \ldots, n\nu_{n-1}$. Hence the relation between $(\nu_1, \ldots, \nu_{n-1})$ and $(\alpha_1, \ldots, \alpha_n)$ is given by
\begin{equation}\label{defnu}
  \nu_j = \frac{1}{n}(\alpha_j - \alpha_{j+1})
\end{equation}
and 
\begin{equation}\label{defalpha}
  \alpha_j = \sum_{i=1}^{n-1} c_{ij} \nu_i, 
  \quad c_{ij} = \begin{cases} n-i, & 1\leq  j \leq i,\\ -i, & i < j \leq n.\end{cases}
\end{equation}
The Plancherel measure is given by (see e.g.\ \cite[p.\ 127]{LM} as well as \cite[Section 4]{LM} which contains a local Weyl law for $GL(n)$)
\begin{displaymath} 
 \prod_{1\leq j < k \leq n} G(\alpha_j - \alpha_k)
\end{displaymath}
which equals \eqref{defspec}. On this occasion we remark that by Stirling's formula 
\begin{displaymath}
    \int_{\| \mu - \nu \| \leq 1}   d_{\text{spec}} \mu \asymp \prod_{1 \leq j \leq k \leq n-1} (1+|\nu_j + \ldots + \nu_k|),
\end{displaymath}
and that the analytic conductor $\mathcal{C}(f \times g)$ of the Rankin-Selberg $L$-function considered   in Theorem 2 satisfies
\begin{displaymath}
  \mathcal{C}(f \times g) \asymp \prod_{1 \leq j \leq k \leq n-1} (1+|\nu_j + \ldots + \nu_k|)^2.
\end{displaymath}
This justifies our earlier remark on the convexity bound implied by Theorem 2.\\


Let $\mathfrak{h}^n$ be the generalized upper half plane as in \cite[p.\ 10]{Go} with coordinates $z = x \cdot y$ where $x \in U_n(\Bbb{R})$, the group of unipotent upper triangular matrices, and $y = \text{diag}(y_1\cdots y_{n-1}, \ldots, y_1, 1)$. It is equipped with a Haar measure 
\begin{displaymath}
   d^{\ast}z = dx \, d^{\ast}{y} = \prod_{i, j} dx_{i, j} \prod_{k=1}^{n-1} y_k^{-k(n-k)} \frac{dy_k}{y_k}. 
\end{displaymath}
For $z \in \mathfrak{h}^n$ put 
\begin{equation}\label{tildez}
  \tilde{z} := w (z^{-1})^t w
\end{equation}
where $w$ is the long Weyl element. 

The Whittaker function $W^{\pm}_{\nu} : \mathfrak{h}^n \rightarrow \Bbb{C}$ is given by (analytic continuation in $\nu = (\nu_1, \ldots, \nu_{n-1})$ of)
\begin{displaymath}
  W_{\nu}^{\pm}(z) = \int_{U_{n}(\Bbb{R})} I_{\nu}(w u z)\overline{\psi}_{\pm}(u) du 
\end{displaymath}
where  $\psi^{\pm}(u) = e(\pm u_{n-1} +  u_{n-2}+ \ldots+ u_1)$ (where $u_{n-1}, \ldots, u_1$ are the entries of the secondary diagonal of $u$) and 
\begin{displaymath}
  I_{\nu}(z) = \prod_{i, j = 1}^{n-1} y_{i}^{b_{ij} (\frac{1}{n}+ \nu_j)}, \quad b_{ij} = \begin{cases} ij, & i+j \leq n,\\ (n-i)(n-j), & i+j \geq n. \end{cases}
\end{displaymath}  
Then we have $W_{\nu}^{\pm}(z) = \psi^{\pm}(x) W_{\nu}(y)$ (the $\pm$ sign at $W_{\nu}(y)$ on the right hand side can be dropped, because the dependence on the sign is only in the first factor). Note that this is \emph{not} the completed Whittaker function, sometimes denoted by $W_{\nu}^{\ast}(z)$ in \cite[Section 5]{Go} and used in \cite{St}. It differs from the completed Whittaker function by a factor
\begin{displaymath}
  \prod_{1 \leq j \leq k \leq n-1} \Gamma_{\Bbb{R}} (1 + n(\nu_j+\ldots + \nu_k)).  
\end{displaymath}
For instance, for $n=2$ we have $W_{\nu}(y) = 2 \pi^{1/2} \cosh(\pi |\nu|/2) \sqrt{y} K_{\nu}(2 \pi y)$, see \cite[p.\ 65]{Go}. The normalizing factor $\cosh(\pi | \nu|/2)$ plays the same role as the factor $\Gamma(k)^{-1/2}$ in \eqref{four}. Whittaker functions for higher rank are not very well understood, but the only information we will need is Stade's formula \cite{St}: for $\nu \in i \Bbb{R}^{n-1}$ define $\alpha \in i\Bbb{R}^n$ as in \eqref{defalpha}, and for $\mu \in i\Bbb{R}^{n-1}$ define $\beta \in i\Bbb{R}^{n}$ correspondingly. Then one has an equality of meromorphic functions in $s$: 
\begin{equation}\label{stade}
\begin{split}
 & \int_{\Bbb{R}_{\geq 0}^{n-1}} W_{\nu}(y) \overline{W_{\mu}}(y) \det(y)^s d^{\ast}y \\
  &= \frac{\prod_{j, k = 1}^n \Gamma_{\Bbb{R}}( s + \alpha_j - \beta_k )}{2  \Gamma_{\Bbb{R}}(ns) \prod_{1 \leq j \leq k \leq n-1} \Gamma_{\Bbb{R}} (1 + n(\nu_j+\ldots + \nu_k))\Gamma_{\Bbb{R}}(1 -n(\mu_j + \ldots + \mu_k)) }. 
  \end{split}
\end{equation}
For real $s \in [1/2, 3/2]$ and  $\mu = \nu + O(1)$ it follows from \eqref{defnu} and \eqref{defalpha} that   the right hand side is  
\begin{equation}\label{simple}
  \asymp  \prod_{1 \leq j \leq k \leq n-1}  \frac{\Gamma_{\Bbb{R}}(s + n(\nu_j+ \ldots + \nu_k))}{\Gamma_{\Bbb{R}}(1+ n(\nu_j+ \ldots + \nu_k))}. 
\end{equation}
In particular for $s=1/2$ we have 
\begin{equation}\label{specgamma}
 \int_{\Bbb{R}_{\geq 0}^{n-1}} W_{\nu}(y) \overline{W_{\mu}}(y) \det(y)^{1/2} d^{\ast}y \asymp \Bigl(\int_{\| \mu - \nu \| \leq 1}   d_{\text{spec}} \mu\Bigr)^{-1/2},  
\end{equation}
cf.\ \eqref{defspec}. \\

Let $f$ be a tempered Hecke-Maa{\ss} form for the group $\Gamma  = SL_n(\Bbb{Z})$ with spectral parameters $\nu = (\nu_1, \ldots, \nu_{n-1}) \in i \Bbb{R}^{n-1}$ and Fourier expansion
\begin{displaymath}
  f(z) = \sum_{\gamma \in U_{n-1}(\Bbb{Z})\backslash SL_{n-1}(\Bbb{Z})} \sum_{m_1=1}^{\infty} \cdots \sum_{m_{n-2} = 1}^{\infty} \sum_{m_{n-1} \not= 0} \frac{A(m_1, \ldots, m_{n-1})}{\prod_{k=1}^{n-1} |m_k|^{k(n-k)/2}} W_{\nu}^{\pm}( M \gamma z)
  \end{displaymath}
  where $\pm = \text{sign} (m_{n-1})$, $M = \text{diag}(m_1\cdots |m_{n-1}|, \ldots, m_1, 1)$, $\gamma$ is embedded in $SL_n(\Bbb{Z})$ as $\left(\begin{smallmatrix} \gamma & \\ & 1 \end{smallmatrix}\right)$ and $A(m_1, \ldots, m_{n-1})$ are Hecke eigenvalues, in particular $A(1, \ldots, 1) = 1$. We have $A(m_1, \ldots, m_{n-1}) = \pm A(m_1, \ldots, -m_{n-1})$. It  follows from  \cite[Theorem 2]{Li} that the Rankin-Selberg $L$-function
  \begin{equation}\label{li}
    L(f \times \bar{f}, s) = \zeta(ns) \sum_{m_1, \ldots, m_{n-1}=1}^{\infty} \frac{|A(m_1, \ldots, m_{n-1})|^2}{(m_1^{n-1} m_2^{n-2} \cdots m_{n-1})^s}
  \end{equation}
 is  bounded by $O(\|\nu\|^{\varepsilon})$ in $\Re s \geq 1 +\varepsilon$, and the same bound holds for its residue at $s=1$.  \\
 
 The space of automorphic forms on $\Gamma \backslash \mathfrak{h}^n$ is equipped with the standard inner product $\langle f, g \rangle = \int_{\Gamma \backslash \mathfrak{h}^n} f(z) \bar{g}(z) d^{\ast}z$.  A fundamental domain for $\Gamma \backslash \mathfrak{h}^n$ is contained in the Siegel set \cite[Prop.\ 1.3.2.]{Go}
 \begin{equation}\label{siegel}
   S := \{z = x \cdot y \in \mathfrak{h}^n \mid 0 \leq x_{ij} \leq 1, y_j \geq \sqrt{3}/2\}. 
 \end{equation}

We introduce the maximal Eisenstein series
\begin{displaymath}
  E^{\ast}(z, s) = \Gamma_{\Bbb{R}}(ns)\zeta(ns) E(z, s)
\end{displaymath}
where
\begin{displaymath}
  E(z, s) :=  \sum_{\gamma \in  P\backslash \Gamma} \det(\gamma z)^s
\end{displaymath}
with
\begin{displaymath}
  P = \left\{ \left(\begin{matrix} * & \cdots & * & *\\ \vdots & \cdots & \vdots &  \vdots\\ * & \cdots & * & *\\ 0 & \cdots & 0 & *\end{matrix}\right) \in SL_n(\Bbb{Z})\right\}. 
\end{displaymath}
It is defined initially for $\Re s > 1$ and can be continued meromorphically to all $s \in \Bbb{C}$. Its Fourier expansion is given explicitly for $n   =2$ and $n=3$ in \cite[p.\ 58, p.\ 226]{Go}. For general $n$ there is an inductive procedure to obtain a certain type of Fourier expansion, see \cite{Te}. We will use this to prove the following bound.  
\begin{lemma}   For $z \in S$ we have
\begin{displaymath}
  E^{\ast}(z, 1/2) \ll \det(z)^{1/2+\varepsilon} + \det (\tilde{z})^{1/2+\varepsilon} 
\end{displaymath}
where $\tilde{z}$ is as in \eqref{tildez}, and the determinant is taken after bringing the matrix back to the canonical Iwasawa form. 
\end{lemma}

\textbf{Proof.} For a positive definite $n \times n$-matrix $M$ and $\Re \rho > n/2$ define the Epstein zeta-function by
\begin{equation}\label{defEp}
  Z(M, \rho) = \frac{1}{2} \sum_{a \in \Bbb{Z}^n \setminus \{0\}} (a^t M a)^{-\rho}. 
\end{equation}
This can be continued meromorphically to all $\rho \in \Bbb{C}$. 
Assume that $M= X^tZX$ where $X \in U_n(\Bbb{R})$ is a unipotent upper triangular matrix and $Z = \text{diag}(z_n, \ldots, z_1)$ is a diagonal matrix with $z_n \gg z_{n-1} \gg \ldots \gg  z_1 \gg 1$. Using the rapid decay of the Bessel $K$-function, it follows by a simple induction from \cite[Theorem 1]{Te} with $n_2 = 1$ and $n_1 = n-1, n-2, \ldots$   that
\begin{equation}\label{Epstein}
 \Gamma_{\Bbb{R}}(2\rho) Z(M, \rho) = \sum_{j=1}^n \Bigl( \Gamma_{\Bbb{R}}(2\rho + 1 - j) \zeta(2\rho+1-j)  + T_j(M, \rho)\Bigr)  (z_1 \cdots z_j)^{-\frac{1}{2}} z_j^{\frac{j}{2}- \rho} . 
\end{equation}
where   $T_j(M, \rho)$ as a function of $\rho$ is holomorphic and  bounded on compact sets, uniformly in $M = X^t Z X$. By \cite[(10.7.3) - (10.7.4)]{Go} we have for $z = x\cdot y \in \mathfrak{h}^n$ the equality\footnote{the factor $\frac{1}{2}$ in \eqref{defEp} is canceled by our definition of $P$ that slightly differs from  the definition in \cite[p.\ 307]{Go}.}
\begin{displaymath}
  E^{\ast}(z, s) =  \det(z)^s \Gamma_{\Bbb{R}}(ns) Z(x^ty^tyx, ns/2).
\end{displaymath}
 Using \eqref{Epstein} with $X = x$, $Z = y^ty = \text{diag}((y_{n-1} \cdots y_1)^2, \ldots, y_1^2, 1)$ and $\rho = ns/2$ and assuming $z = x \cdot y \in S$ we find 
\begin{displaymath}
  E^{\ast}(z, s) = \sum_{j=1}^n \Bigl(\Gamma_{\Bbb{R}}(ns + 1 - j) \zeta(ns+1-j)   + \tilde{T}_j(z, s)\Bigr) \prod_{i=1}^{n-1} y_i^{a_{ij}(s)}
\end{displaymath} 
where
\begin{displaymath}
   a_{ij}(s) = \begin{cases} (n-i)s, & 1 \leq j \leq i,\\ i(1-s), & i < j \leq n, \end{cases}
\end{displaymath}
and $\tilde{T}_j(z, s)$ is holomorphic and uniformly bounded for $z \in S$. Specializing to $s=1/2$ and taking residues if necessary, we conclude
\begin{displaymath}
  E^{\ast}(z, 1/2) \ll \sum_{j=1}^n   \prod_{i=1}^{n-1} y_i^{a_{ij}(1/2) +\varepsilon} \ll \det(z)^{1/2+\varepsilon} + \det (\tilde{z})^{1/2+\varepsilon} 
\end{displaymath}
for $z \in S$. \qed\\

Again we can use the Rankin-Selberg unfolding method to compute the norm of a Maa{\ss} form $f$, see \cite[Section 12.1]{Go}. The Eisenstein series $E(z, s)$ has a pole at $s=1$ with constant residue \cite[p.\ 483]{Te}, hence
\begin{displaymath}
  \| f \|^2 \asymp  \underset{s=1}{\text{res}} \int_{\Gamma\backslash \mathfrak{h}^n} |f(z)|^2 E(z, s) d^{\ast}z \asymp   \underset{s=1}{\text{res}} L(f \times \bar{f}, s)  \int_{\Bbb{R}_{\geq 0}^{n-1}} W_{\nu}(y) \overline{W_{\mu}}(y) \det(y)^s d^{\ast}y.
\end{displaymath}
We can compute the integral by Stade's formula \eqref{stade} at $s=1$. By \eqref{simple} and recalling Li's bound for \eqref{li}, we see
\begin{equation}\label{sizenorm}
  \| f \|^2 \asymp  \underset{s=1}{\text{res}} L(f \times \bar{f}, s)  \ll \| \nu \|^{\varepsilon}. 
\end{equation}
This shows the usefulness of our normalization of Whittaker functions. 
  
 \section{Proofs of Theorems 2 and 3} 
 
We are now ready to imitate the proof from Section 2 in the general case. Let $f$ be as in Theorem 2, and use the notation from Theorem 2. Then by Bessel's inequality we have 
\begin{equation}\label{bessel}
 \| f E^{\ast}(., s) \|^2 \geq \sum_{\| \mu - \nu\| \leq 1} \frac{1}{\| g\|^2} |\langle f E^{\ast}(., s), g\rangle|^2.
\end{equation}
Unfolding as above we obtain (first in $\Re s > 1$)
\begin{displaymath}
  \langle f E^{\ast}(., s), g\rangle = \Gamma_{\Bbb{R}}(ns) L(f \times \bar{g}, s) \int_{\Bbb{R}_{\geq 0}^{n-1}} W_{\nu}(y) \overline{W_{\mu}}(y) \det(y)^s d^{\ast}y.
\end{displaymath}
We use  Stade's formula\footnote{It is not hard to see that  the Whittaker integral  is negligible unless $\mu \approx \nu$, so that the inequality \eqref{bessel} does not lose much.} \eqref{stade},  extend both sides meromorphically to all $s \in \Bbb{C}$ and   specialize $s=1/2$.  In the range $\|\mu - \nu \| \leq 1$, we use \eqref{sizenorm} for $g$ and   \eqref{specgamma} to conclude 
\begin{displaymath}
\begin{split}
   \sum_{\| \mu - \nu \| \leq 1} |L(f \times \bar{g}, 1/2)|^2& \ll \|\nu\|^{\varepsilon} \sum_{\| \mu - \nu \| \leq 1}\frac{1}{\|g\|^2}  |L(f \times \bar{g}, 1/2)|^2\\
& \ll  \Bigl(\int_{\| \mu - \nu \| \leq 1}   d_{\text{spec}} \mu\Bigr)^{1+\varepsilon} \| f E^{\ast}(., 1/2)\|^2.
\end{split}
\end{displaymath}
Let $\mathcal{F}$ be a fundmental domain for $\Gamma \backslash \mathfrak{h}^n$ contained in the Siegel set $S$ defined in \eqref{siegel}. Let $\tilde{f}(z) = f(w (z^{-1})^t w)$ be the dual Maass form. Then by Lemma 1,
\begin{displaymath}
  \| f E^{\ast}(., 1/2) \|^2 \ll  \int_{\mathcal{F}} (|f(z)|^2 + |\tilde{f}(z)|^2) \det(z)^{1+\varepsilon} d^{\ast}z 
   \leq \int_{\mathcal{F}} (|f(z)|^2 + |\tilde{f}(z)|^2) E(z, 1+\varepsilon) d^{\ast}z. 
 \end{displaymath}   
    The first term equals
    \begin{displaymath}
     \sum_{m_1, \ldots, m_{n-1}=1}^{\infty} \frac{|A(m_1, \ldots, m_{n-1})|^2}{(m_1^{n-1} m_2^{n-2} \cdots m_{n-1})^{1+\varepsilon}} \int_{\Bbb{R}_{\geq 0}^{n-1}} W_{\nu}(y) \overline{W_{\nu}}(y) \det(y)^{1+\varepsilon} d^{\ast}y,
\end{displaymath}  
and the second term is similar with indices interchanged. Appealing to Stade's formula one last time, and using the uniform bound for the $L$-series, the preceding quantity is $O( \| \nu\|^{\varepsilon})$ and the proof is   complete. \\

The proof of Theorem 3 is the same and uses the deep connection between triple products and central values of $L$-functions, as developed in \cite{Wa, Ic, Wo}. Fix a Hecke-Maa{\ss} cusp form $h$ with spectral parameter $\lambda \in i\Bbb{R}$, and let $f, g$ be as before. We have by \eqref{sizenorm}
\begin{displaymath}
 \frac{1}{|\nu|^{1+\varepsilon}} \sum_{|\mu - \nu| \leq 1} L(f \times \bar{g} \times h, 1/2)  \ll_h  \sum_{|\mu - \nu| \leq 1} \frac{1}{\| f \|^2 \|g \|^2 \| h \|^2} \Lambda^{\ast}(f \times \bar{g} \times h, 1/2).
 \end{displaymath}
 Here 
 \begin{displaymath}
   \Lambda^{\ast}(f \times \bar{g} \times h, s) = \Lambda(f \times \bar{g} \times h, s) \cosh\left(\frac{\pi |\nu|}{2}\right)^2 \cosh\left(\frac{\pi |\mu|}{2}\right)^2
\cosh\left(\frac{\pi |\lambda|}{2}\right)^2
 \end{displaymath}
 is the ``normalized" completed $L$-function (so that $\|f \|, \|g\|, \|h \| \approx 1$). 
 By Watson's formula, this is
 \begin{displaymath}
   \sum_{|\mu - \nu | \leq 1} |\langle fh, g\rangle|^2 \leq \| f h \|^2 \leq \| h \|_{\infty} \| f \|^2  \ll_h |\nu |^{\varepsilon}. 
\end{displaymath}

\end{document}